\documentclass{icmart}

\contact[verbit@verbit.ru]{7 Vavilova Str., Moscow, Russia, 117312}

\newtheorem{theorem}{Theorem}[section]
\newcommand{\Theorem}{\begin{theorem}}
\newcommand{\et}{\end{theorem}}

\newtheorem{proposition}[theorem]{Proposition}
\newcommand{\Proposition}{\begin{proposition}}
\newcommand{\ep}{\end{proposition}}

\newtheorem{conjecture}[theorem]{Conjecture}
\newcommand{\Conjecture}{\begin{conjecture}}
\newcommand{\econ}{\end{conjecture}}

\newtheorem{corollary}[theorem]{Corollary}
\newcommand{\Corollary}{\begin{corollary}}
\newcommand{\eco}{\end{corollary}}

\newtheorem{claim}[theorem]{Claim}
\newcommand{\Claim}{\begin{claim}}
\newcommand{\ec}{\end{claim}}

\theoremstyle{definition}
\newtheorem{definition}[theorem]{Definition}
\newcommand{\Definition}{\begin{definition}}
\newcommand{\ed}{\end{definition}}

\newtheorem{remark}[theorem]{Remark}
\newcommand{\Remark}{\begin{remark}}
\newcommand{\er}{\end{remark}}

\newtheorem{example}[theorem]{Example}
\newcommand{\Example}{\begin{example}}
\newcommand{\ex}{\end{example}}

\newtheorem{question}[theorem]{Question}
\newcommand{\Question}{\begin{question}}
\newcommand{\eq}{\end{question}}

\renewcommand{\phi}{\varphi}
\renewcommand{\epsilon}{\varepsilon}
\renewcommand{\geq}{\geqslant}
\renewcommand{\leq}{\leqslant}
\renewcommand{\Re}{\operatorname{Re}}
\renewcommand{\Im}{\operatorname{Im}}

\newcommand{\Z}{{\Bbb Z}}
\newcommand{\C}{{\Bbb C}}
\newcommand{\R}{{\Bbb R}}
\newcommand{\Q}{{\Bbb Q}}
\newcommand{\calo}{{\cal O}}
\newcommand{\1}{{\sqrt 1}}

\newcommand{\End}{\operatorname{End}}
\newcommand{\St}{\operatorname{St}}
\newcommand{\Sym}{\operatorname{Sym}}

\newcommand{\Id}{\operatorname{Id}}
\newcommand{\Vol}{\operatorname{Vol}}

\newcommand{\Aut}{\operatorname{Aut}}

\newcommand{\Lie}{\operatorname{Lie}}
\newcommand{\Kah}{\operatorname{Kah}}
\newcommand{\diam}{\operatorname{\sf diam}}
\renewcommand{\max}{{\operatorname{\sf max}}}
\newcommand{\Diff}{\operatorname{Diff}}
\newcommand{\Teich}{\operatorname{Teich}}
\newcommand{\Comp}{\operatorname{Comp}}
\newcommand{\Per}{\operatorname{\sf Per}}
\newcommand{\Perspace}{\operatorname{{\Bbb P}\sf er}}
\newcommand{\arrow}{{\:\longrightarrow\:}}

\newcommand{\goth}{\mathfrak}
\newcommand{\restrict}[1]{{\left|_{{#1}}\right.}}
\def\blacksquare{\hbox{\vrule width 5pt height 5pt depth 0pt}}
\def\endproof{\blacksquare}

\def\eqref#1{(\ref{#1})}
\renewcommand{\theequation}{\thesection.\arabic{equation}}

\makeatletter
\renewcommand\subsection{\@startsection{subsection}{2}{\z@}%
                                     {12pt\@plus -1ex \@minus -.2ex}%
                                     {.7ex \@plus .2ex}%
                                     {\normalfont\large\bfseries\boldmath}}
\makeatother

\title[Teichm\"uller spaces, ergodic theory and global Torelli theorem]%
{Teichm\"uller spaces, ergodic theory \\and global Torelli theorem}

\author[Misha Verbitsky]{Misha Verbitsky}

\begin{document}

\begin{abstract}
A Teichm\"uller space $\Teich$ is a quotient of the space of 
all complex structures on a given manifold $M$
by the connected components of the group of 
diffeomorphisms. The mapping class group $\Gamma$ of $M$ 
is the group of connected components of the
diffeomorphism group. The moduli problems
can be understood as statements about the 
$\Gamma$-action on $\Teich$.
I will describe the mapping class group and
the Teichm\"uller space for a hyperk\"ahler manifold.
It turns out that this action is ergodic. We use
the ergodicity to show that a hyperk\"ahler manifold
is never Kobayashi hyperbolic.
\end{abstract}



\maketitle

\tableofcontents

\hfill

This talk is based on two papers,
\cite{_V:Torelli_} and \cite{_V:ergodic_}. In these papers one can find 
details, examples, and rigorous proofs omitted here.


\section{Teichm\"uller spaces}


\subsection{Teichm\"uller spaces and period maps}

The notion of Teichm\"uller spaces has a long history
since its discovery by Teichm\"uller in 1944 
(\cite{_Teichmuller:1944_}) and further development 
by Ahlfors, Bers and others. However, it is rarely applied
to complex manifolds of dimension $>1$. It turns out
that this notion is interesting and useful for 
many purposes of complex geometry in any dimension.

\begin{definition}
Let $M$ be a smooth manifold. 
An {\bf almost complex structure} is an operator
$I:\; TM \arrow TM$ which satisfies $I^2 = - \Id_{TM}$.
An almost complex structure is {\bf  integrable}
if $\forall X,Y \in T^{1,0}M$, one has $[X,Y]\in T^{1,0}M$.
In this case $I$ is called {\bf  a complex structure operator}.
A manifold with an integrable almost complex structure
is called {\bf  a complex manifold}. 
\end{definition}

\begin{definition}
{\bf The space of almost complex structures}
is an infinite-dimensional Fr\'echet manifold $X_M$  of all tensors
$I\in \End(TM)$ satisfying $I^2 = - \Id_{TM}$. 
Similarly, one considers {\bf the group
of diffeomorphisms} as a Fr\'echet Lie group. 
\end{definition}

\begin{remark} Definition of Fr\'echet manifolds and Fr\'echet spaces
and many results on the geometry of infinite-dimensional
manifolds can be found in \cite{_Hamilton:inverse_}.
\end{remark}

\begin{definition}
Let $M$ be a compact complex manifold, and 
$\Diff_0(M)$ a connected component of its diffeomorphism group
({\bf the group of isotopies}). Denote by $\Comp$
the space of complex structures on $M$, considered with the
topology induced from a Fr\'echet manifold of almost
complex structures, and let $\Teich:=\Comp/\Diff_0(M)$ 
be the quotient space with the quotient topology. We call 
it {\bf the Teichm\"uller space.}
\end{definition}

\begin{remark} 
When the complex manifold $M$ admits a certain geometric structure,
such as K\"ahler, or hyperk\"ahler structure, it is natural
to consider the Teichm\"uller space of complex structure 
compatible with (say) K\"ahler structure. Consider the 
open subset $\Comp_K\subset \Comp$ of all complex structures
$I$ such that $(M,I)$ admits a K\"ahler structure. The corresponding
Teichm\"uller space is $\Teich_K:=\Comp_K/\Diff_0$. When
working with the Teichm\"uller space of hyperk\"ahler manifolds,
or a torus, we shall always restrict ourselves to $\Comp_K$
and $\Teich_K$.
\end{remark}

Results of Kuranishi about local structure of deformation
spaces can be summarized as a statement 
about local structure of $\Comp$ as follows
(\cite{_Kuranishi:new_,_Kuranishi:note_,_Douady:Bourbaki_}).

\begin{theorem}
Let $M$ be a compact complex manifold and $I\in \Comp$.
Then there exists an open neighbourhood $U\ni I$ in $\Comp$
and a neighbourhood $R$ of unit in $\Diff_0$ satisfying the
following. Consider the quotient $U/R$ of $U$ by an equivalence
relation generated by $x\sim wy$, for all $x,y\in U$ and $w\in R$.
Then $U/R$ is a complex variety, equipped with a natural holomorphic
embedding to $H^1(TM)$.
\end{theorem}

\begin{remark}\label{_Kuranishi_Remark_}
The quotient space  $U/R$ obtained by Kuranishi is
called {\bf the Kuranishi space}. Let $\Teich(U)$ be an
image of $U$ in the Teichm\"uller space. Clearly, the
Kuranishi space admits a surjective, continuous 
map to $\Teich(U)$. It is not entirely clear whether
this map is always a homeomorphism. However, if it
is always a homeomorphism, for a given $M$, the space
$\Teich$ acquires a structure of a complex variety.
\end{remark}

As shown by F. Catanese (\cite[Proposition 15]{_Catanese:moduli_}),
for K\"ahler manifolds with trivial canonical bundle,
e.g. for the hyperk\"ahler manifolds,
the Teichm\"uller space is locally isomorphic
to the Kuranishi moduli space, hence it is a complex
variety. In this case it is actually a complex manifold, by
Bogomolov-Tian-Todorov theorem (\cite{_Bogomolov_81_,_Tian_,_Todorov_}).

It is not clear if this is true for
a general complex manifold; in the present
work we deal with hyperk\"ahler manifolds,
which are Calabi-Yau.

\begin{question}
Consider a compact complex manifold $M$, and
let $\Teich$ be its Teichm\"uller space.
Can we equip $\Teich$ with a structure of a
complex variety (possibly non-Hausdorff), in 
a way which is compatible with the local charts
obtained from the Kuranishi theorem?
\end{question}

When $M$ is a torus, or a hyperk\"ahler manifold,
$\Teich$ is a complex manifold which can be described 
explicitly (Theorem \ref{_torus_Torelli_Theorem_},
Theorem \ref{_Glo_Tor_Theorem_}).
However, even for a hyperk\"ahler manifold,
$\Teich$ is not Hausdorff.

\begin{claim}
Assume that $M$ is K\"ahler, and $\Teich$ the Teichm\"uller
space of all complex structures of K\"ahler type on $M$. 
For a given $I\in \Teich$, choose a representative
$\tilde I\in \Comp$. Then the Hodge decomposition
$H^*(M)=\bigoplus H^{p,q}(M,\tilde I)$ is independent from the
choice of $\tilde I$.
\end{claim}

{\bf Proof:} The ambiguity of a choice of $\tilde I$ lies
in $\Diff_0$. However, $\Diff_0$ acts trivially on $H^*(M)$.
\endproof

\hfill

This elementary claim allows one to define
{\bf the period map}.

\begin{definition}\label{_period_Definition_}
Lety $M$ be a K\"ahler manifold, $\Teich$ 
its Teichm\"uller space, and $\Per$ the map
associating to $I$ the Hodge decomposition
$H^*(M)=\bigoplus H^{p,q}(M,I)$. Then
$\Per$ is called {\bf the period map} of $M$.
\end{definition}

\begin{remark}\label{_Universal_fi_Remark_}
Consider the product $\Comp\times M$ trivially fibered over
$\Comp$. The fibers of $\pi:\; M \times \Comp\arrow \Comp$
can be considered as complex manifolds, with complex
structure at $I\in \Comp$ given by $I$. This complex 
structure is clearly $\Diff_0$-invariant, giving a complex
structure on the fibers of the quotient fibration
$(M\times \Comp)/\Diff_0 \arrow \Teich$. At each $I\in \Teich$,
the fiber of this fibration (called {\bf the universal fibration})
is isomorphic to $(M,I)$. 
\end{remark}

\subsection{Marked moduli spaces}

A more conventional approach to the moduli problem
goes as follows. Given a complex manifold $M$, 
one defines {\bf the deformation functor}
from marked complex spaces to sets as a functor
mapping a complex space $(B,x)$ to the set of equivalence
classes of deformations $\pi:\; {\cal X}\arrow B$ of $M$ over $B$
with $M$ identified with the fiber of $\pi$ at $x$.

If the deformation functor is representable by 
a complex space, this space is called {\bf the fine moduli space}
of deformations of $M$. 

Usually, the fine moduli space
does not exist. In this case, one considers the category of
natural transformations from the deformation functor
to representable functors. The initial object in this category
is called {\bf the coarse moduli space}. The points of coarse
moduli are identified with equivalence classes of deformations
of $M$.

In this setup, an analogue of Teichm\"uller space can be defined
as follows. Fix an abelian group which is isomorphic to
$H^*(M,\Z)$, and define {\bf a marked manifold} as a pair 
$(M, \phi:\; V \tilde \arrow H^*(M,\Z))$, where $M$ is a complex
manifold, and $\phi$ a group isomorphism. In the same way as above, 
one defines a coarse moduli space of deformations
of marked manifolds. To compare this space with Teichm\"uller space,
consider a subgroup group $\Gamma_0$ of mapping class group
which acts trivially on cohomology. Clearly, the points of 
$\Teich/\Gamma_0$ are in bijective correspondence with the
equivalence classes of marked complex structures on $M$.

Given a coarse marked moduli space $W$,
one obtains the tautological map $W\arrow \Teich/\Gamma_0$, 
by construction continuous. For hyperk\"ahler manifolds
(or compact tori), this map is a diffeomorphism on each
connected component (\cite[Corollary 4.31]{_V:Torelli_}).


\section{Torelli theorem}


\subsection{Torelli theorem: an introduction}

Torelli theorems are a broad class of results which
describe the Teichm\"uller spaces in terms of
the period maps (Definition \ref{_period_Definition_}).

The name originates with Ruggiero Torelli, who 
has shown that it is possible to reconstruct 
a Riemann surface from its Jacobian (\cite{_Torelli_}). 
The term ``Torelli theorems'' is due to Andr\'e Weil
(\cite{_Weil_}), who gave a modern proof of this classical
result, and explained its possible generalizations.

One may distinguish between the ``local Torelli theorem'',
where a local structure of deformation space is described
in terms of periods, and ``global Torelli'', where the
Teichm\"uller space is described globally.

Weil, who was the first to define and 
study K3 surfaces, spent much time trying 
to prove the Torelli theorem for K3 surfaces,
but it was notoriously difficult. Its local version
is due to Tjurina, Piatetski-Shapiro and Shafarevich 
(\cite{_Tjurina_,_Piatetski_Shapiro_Shafarevich_}).
The local Torelli was generalized by Bogomolov to
hyperk\"ahler manifolds (\cite{_Bogomolov_78_})
and by Bogomolov-Tian-Todorov
to Calabi-Yau manifolds (\cite{_Bogomolov_81_,_Tian_,_Todorov_}), 
building foundation for the theory of Mirror Symmetry.

In dimension $>1$, the global Torelli theorem 
was known only for compact tori (where it is essentially trivial)
and the K3 surfaces, where it was proven by Kulikov in 1977
(\cite{_Kulikov_}), and then improved many times during the 1980-ies
(\cite{_Todorov:K3_,_Beauville:K3_,_Loojenga_,_Siu:K3_,_Friedman_}).

\subsection{Birational Teichm\"uller space}

In what follows, {\bf a hyperk\"ahler manifold}
is a compact complex manifold admitting a K\"ahler structure
and a holomorphic symplectic form.

In generalizing global Torelli to more general hyperk\"ahler manifolds,
two problems were apparent. First of all, bimeromorphic hyperk\"ahler
manifolds have the same periods, hence the period map cannot
distinguish between them. However, for $\dim_\C M>0$,
birational holomorphically symplectic manifolds can be non-isomorphic 
(\cite{_Debarre:Torelli_}). 

Another (mostly psychological) difficulty is based on 
attachment to moduli spaces, as opposed to marked moduli
or Teichm\"uller spaces. For a K3, one can reconstruct
a K3 from its Hodge structure, and this gives an identification
between the moduli and the space of Hodge structures. In 
bigger dimension, one has to use the Teichm\"uller space.
Indeed, for some classes of hyperk\"ahler manifolds,
the group $O(H^2(M,\Z)$ of Hodge isometries of cohomology 
is strictly bigger than the image of the mapping class group.
This gives elements $\gamma\in O(H^2(M,\Z)$ acting non-trivially
on the Teichm\"uller space in such a way that the complex manifolds 
$(M,I)$ and $(M, \gamma(I))$ are not birationally equivalent
(\cite{_Namikawa:Torelli_,_Markman:constra_}). However, their
Hodge structures {\em are} equivalent, by construction.
This example explains the necessity of using the 
Teichm\"uller spaces (or marked moduli) 
to state the Torelli theorem: its Hodge-theoretic 
version is often false.

For Teichm\"uller spaces, the Torelli theorem is a statement
about the period map (Definition \ref{_period_Definition_}). 
Ideally, we want the period map to give a diffeomporphism 
between $\Teich$ and the corresponding space of Hodve structures.
This is what happens for a compact torus.

\begin{theorem}\label{_torus_Torelli_Theorem_}
Let $M$ be a compact torus, $\dim_\R M=2n$, and $\Teich$ the
Teichm\"uller space of all complex structures of K\"ahler type on
$M$. Denote by $\Perspace$ the space $SL(2n, \R)/SL(n,\C)$
of all Hodge structures of weight one on $H^1(M,\C)$,
that is, the space of all complex operators on $H^1(M,\R)$
compatible with the orientation. Then the period map
$\Per:\; \Teich\arrow \Perspace$ is a diffeomorphism
on each connected component of $\Perspace$. \endproof
\end{theorem}

Unfortunately, this ideal situation is almost never realized. 
Even in the simplest cases (such as for hyperk\"ahler manifolds),
the Teichm\"uller space is no longer Hausdorff. 
However, in some situations it is still possible to 
deal with non-Hausdorff points.

\Remark {\bf A non-Hausdorff manifold} is a topological
space locally diffeomorphic to $\R^n$ (but not necessarily Hausdorff).
\er

\Definition Let $X$ be a topological space, and $X\stackrel \phi \arrow X_0$ 
a continuous surjection. The space $X_0$ is called 
{\bf a Hausdorff reduction} of $X$ if any
continuous map $X\arrow X'$ to a Hausdorff space
is factorized through $\phi$.
\ed

\Definition
Let $M$ be a topological space. We say that $x, y \in M$
are {\bf non-separable} (denoted by $x\sim y$)
if for any open sets $V\ni x, U\ni y$, one has $U \cap V\neq \emptyset$.
\ed

\Remark
Suppose that $\sim$ is an equivalence relation, and
the quotient $M/\sim$ is Hausdorff. Then $M/\sim$ is a 
Hausdorff reduction of $M$.
\er

Unfortunately, this notion cannot be applied 
universally. Firstly, $\sim$ is not always an equivalence relation;
and secondly, even if $\sim$ is equivalence, the $M/\sim$ 
is not always Hausdorff.
Fortunately, for Teichm\"uller space of a hyperk\"ahler manifold,
Hausdorff reduction can be defined, using the following
theorem due to D. Huybrechts (\cite{_Huybrechts:cone_}).

\Theorem
If $I_1$, $I_2\in \Teich$ are non-separate points, 
then $(M, I_1)$ is birationally equivalent to $(M, I_2)$.
\et

Using this result and geometry of the period map
(Bogomolov's local Torelli theorem), it is elementary
to show that the quotient $\Teich_b:=\Teich/\sim$ is a Hausdorff manifold.
This quotient is called {\bf the birational Teichm\"uller space}
of a hyperk\"ahler manifold. 

Global Torelli theorem implies that 
for hyperk\"ahler manifolds the period map induces a diffeomorphism
between the Hausdorff reduction of the Teichm\"uller space and
the appropriate period domain.


\section{Hyperk\"ahler manifolds and Bogomolov-Beauville-Fujiki form}


\subsection{Hyperk\"ahler manifolds: definition and examples}

The standard definition of hyperk\"ahler manifolds is rather differential
geometric. It is, indeed, synonymous with ``holomorphic symplectic'',
but this synonymity follows from Calabi-Yau theorem.
For more details about hyperk\"ahler manifolds, please see
\cite{_Beauville:hk_} or \cite{_Besse:Einst_Manifo_}.

\Definition
A {\bf hyperk\"ahler structure} on a manifold $M$
is a Riemannian structure $g$ and a triple of complex
structures $I,J,K$, satisfying quaternionic relations
$I\circ J = - J \circ I =K$, such that $g$ is K\"ahler
for $I,J,K$.
\ed

\Remark
This is equivalent to $\nabla I=\nabla J = \nabla K=0$:
the parallel translation along the connection preserves $I, J,K$.
\er

\Remark A hyperk\"ahler manifold has three symplectic forms:
$\omega_I:=  g(I\cdot, \cdot)$, $\omega_J:=  g(J\cdot, \cdot)$,
$\omega_K:=  g(K\cdot, \cdot)$.
\er

\Definition
Let $M$ be a Riemannian manifold, $x\in M$ a point.
The subgroup of $GL(T_xM)$ generated by parallel 
translations (along all paths) is called {\bf 
the holonomy group} of $M$.
\ed

\Remark
A hyperk\"ahler manifold can be defined as a manifold which has
holonomy in $Sp(n)$ (the group of all endomorphisms preserving
$I,J,K$).
\er

\Definition
A holomorphically symplectic manifold 
is a complex manifold equipped with non-degenerate, holomorphic
$(2,0)$-form.
\ed

\Remark
 Hyperk\"ahler manifolds are holomorphically symplectic.
Indeed, $\Omega:=\omega_J+\1\omega_K$ is a holomorphic symplectic
form on $(M,I)$.
\er

\Theorem[(Calabi-Yau)] 
A compact, K\"ahler, holomorphically symplectic manifold
admits a unique hyperk\"ahler metric in any K\"ahler class.
\et

\Remark
 For the rest of this talk, a hyperk\"ahler manifold
means a compact complex manifold admitting a K\"ahler
structure and a holomorphically symplectic structure.
\er

\Definition
A hyperk\"ahler manifold $M$ is called
{\bf simple}, or {\bf IHS}, if $\pi_1(M)=0$, $H^{2,0}(M)=\C$.
\ed

The rationale for this terminology comes from Bogomolov's
decomposition theorem.

\Theorem[(Bogomolov, \cite{_Bogomolov:decompo_})]
 Any hyperk\"ahler manifold admits a finite covering
which is a product of a torus and several 
simple hyperk\"ahler manifolds.
\et 

 Further on, all hyperk\"ahler manifolds
are assumed to be simple.

\Remark\label{_simple_holo_Remark_}
A hyperk\"ahler manifold is simple if and only if its holonomy
group is $Sp(n)$, and not a proper subgroup of $Sp(n)$ 
(\cite{_Besse:Einst_Manifo_}).
\er

\Example Take a 2-dimensional complex torus $T$,
then the singular locus of $T/{\pm1}$ is 16 points locally 
of form $\C^2/{\pm1}$.
Its resolution by blow-up  is called 
{\bf a Kummer surface}. It is not hard to see that 
it is holomorphically symplectic.
\ex

\Definition
A K3 surface is a hyperk\"ahler manifold which is
diffeomorphic to a Kummer surface.
\ed

In real dimension 4, the only compact hyperk\"ahler manifolds
are tori and K3 surfaces, as follows from the 
Kodaira-Enriques classification.

\Definition
A {\bf Hilbert scheme} $M^{[n]}$ of a complex surface $M$ is
a classifying space of all ideal sheaves $I\subset \calo_M$ 
for which the quotient $\calo_M/I$ has dimension $n$
over $\C$.
\ed

\Remark 
A Hilbert scheme is obtained as a resolution of singularities
of the symmetric power $\Sym^n M$.
\er

\Theorem[(Fujiki, Beauville)] A Hilbert scheme of
a hyperk\"ahler manifold of real dimension 2 is hyperk\"ahler.
\et

\Example
Let $T$ be a torus. Then $T$ acts on its Hilbert scheme
freely and properly by translations. For $n=2$, the quotient $T^{[n]}/T$
is a Kummer K3-surface. For $n>2$, a universal covering
of $T^{[n]}/T$ is called {\bf a generalized Kummer variety}. 
\ex

\Remark There are 2 more ``sporadic'' examples
of compact hyperk\"ahler manifolds, constructed by K. O'Grady
(\cite{_O_Grady_}).
All known simple hyperkaehler manifolds are
these 2 and the two series: Hilbert schemes of K3 and
generalized Kummer.
\er

\subsection{Bogomolov-Beauville-Fujiki form and the mapping class group}

\Theorem
(Fujiki, \cite{_Fujiki:HK_})\\
Let $\eta\in H^2(M)$, and $\dim M=2n$, where $M$ is
hyperk\"ahler. Then $\int_M \eta^{2n}=c q(\eta,\eta)^n$,
for some primitive integer quadratic form $q$ on $H^2(M,\Z)$,
and $c>0$ a positive rational number, called {\bf Fujiki constant}.
\et

\Definition
This form is called
{\bf  Bogomolov-Beauville-Fujiki form}.  It is defined
by the Fujiki's relation uniquely, up to a sign. The sign is determined
from the following formula (Bogomolov, Beauville)
\begin{align*}  \lambda q(\eta,\eta) &=
   \int_X \eta\wedge\eta  \wedge \Omega^{n-1}
   \wedge \bar \Omega^{n-1} -\\
 &-\frac {n-1}{n}\left(\int_X \eta \wedge \Omega^{n-1}\wedge \bar
   \Omega^{n}\right) \left(\int_X \eta \wedge \Omega^{n}\wedge \bar \Omega^{n-1}\right)
\end{align*}
where $\Omega$ is the holomorphic symplectic form, and 
$\lambda>0$.
\ed

\Remark
The BBF form $q$ has signature $(b_2-3,3)$.
It is negative definite on primitive forms, and positive
definite on $\langle \Omega, \bar \Omega, \omega\rangle$,
 where $\omega$ is a K\"ahler form. 
\er

Using the BBF form, it is possible to describe the automorphism
group of cohomology in a very convenient way.

\Theorem
Let $M$ be a simple hyperk\"ahler manifold, and 
$G\subset GL(H^*(M))$ a group of automorphisms of its cohomology
algebra preserving the Pontryagin classes.
Then $G$ acts on $H^2(M)$ preserving the BBF form. Moreover,
the map $G\arrow O(H^2(M, \R), q)$ is surjective on a connected
component, and has compact kernel.
\et

{\bf Proof. Step 1:} Fujiki formula 
$v^{2n}= q(v,v)^n$ implies that
$\Gamma_0$ preserves the Bogomolov-Beauville-Fujiki
up to a sign.  The sign is fixed, if $n$ is odd.

{\bf Step 2:} For even $n$, the sign is also fixed. 
Indeed, $G$ preserves $p_1(M)$, and (as Fujiki has shown
in \cite{_Fujiki:HK_}),
$v^{2n-2}\wedge p_1(M)= q(v,v)^{n-1} c$, 
for some $c\in \R$. The constant $c$ is positive, 
because the degree of $c_2(B)$ is positive
for any non-trivial Yang-Mills bundle with $c_1(B)=0$.

{\bf  Step 3:} ${\goth o}(H^2(M, \R), q)$
acts on $H^*(M, \R)$ by derivations preserving 
Pontryagin classes (\cite{_Verbitsky:coho_announce_}). Therefore 
$\Lie(G)$ surjects to ${\goth o}(H^2(M, \R), q)$.

{\bf  Step 4:} The kernel $K$ of the map
$G \arrow G\restrict{H^2(M,\R)}$ is compact,
because it commutes with the Hodge decomposition and
 Lefschetz ${\goth sl}(2)$-action, hence preserves
the Riemann-Hodge form, which is positive definite.
\endproof 

\hfill

Using this result, the mapping class group can also be computed.
We use a theorem of D. Sullivan, who expressed the mapping group
in terms of the rational homotopy theory, and
expressed the rational homotopy in terms of the 
algebraic structure of the de Rham algebra.

\Theorem (Sullivan, %
\cite[Theorem 10.3, Theorem 12.1, Theorem 13.3]{_Sullivan:infinite_})\\
Let $M$ be a compact, simply connected 
K\"ahler manifold, $\dim_\C M\geq 3$. Denote by $\Gamma_0$ the group
of automorphisms of an algebra $H^*(M, \Z)$
preserving the Pontryagin classes $p_i(M)$. 
Then the natural map 
$\Diff(M)/\Diff_0\arrow \Gamma_0$ has finite kernel,
and its image has finite index in $\Gamma_0$.
\et

As a corollary of this theorem, we obtain a similar result about
hyperk\"ahler manifolds.

\Theorem
Let $M$ be a simple hyperk\"ahler manifold,
and $\Gamma_0$ the group
of automorphisms of an algebra $H^*(M, \Z)$
preserving the Pontryagin classes $p_i(M)$.  Then 
\begin{description}
\item[(i)]  $\Gamma_0\restrict{H^2(M,\Z)}$ is a finite index 
subgroup of $O(H^2(M, \Z), q)$.
\item[(ii)] The map $\Gamma_0\arrow O(H^2(M, \Z), q)$
 has finite kernel. 
\end{description}\endproof
\et

We obtained that the mapping group is {\bf arithmetic} (commensurable
to a subgroup of integer points in a rational Lie group).

\hfill

As follows from \cite[Theorem 2.1]{_Huybrechts:finiteness_}, 
there are only finitely many connected components  of $\Teich$.
 Let $\Gamma^I$ be the group of elements of mapping class 
group preserving a connected component of Teichm\"uller 
space containing $I\in \Teich$. Then $\Gamma^I$
is also arithmetic. Indeed, it has finite index in $\Gamma$.

\Definition\label{_Monodro_Definition_}
The image of $\Gamma^I$ in $GL(H^2(M,\Z))$ is called {\bf monodromy group of
a manifold}. 
\ed

\Remark
The monodromy group can also be obtained as a group generated
by monodromy of all Gauss-Manin local system for all
deformations of $M$ (\cite[Theorem 7.2]{_V:Torelli_}). 
This explains the term. This notion was defined
and computed in many special cases by E. Markman 
(\cite{_Markman:mono_}, \cite{_Markman:constra_}).
\er


\section{Global Torelli theorem}


\subsection{Period map}

To study the moduli problem, one should understand
the mapping class group (described above) and the Teichm\"uller
space. It turns out that the birational Teichm\"uller space
has a very simple description in terms of the period map,
inducing a diffeomorphism $\Teich_b\arrow 
\frac{SO(b_2-3,3)}{SO(b_2-3,1)\times SO(2)}$
on each connected component of $\Teich_b$.

\Definition Let 
$\Per:\; \Teich \arrow {\Bbb P}H^2(M, \C)$
map $J$ to a line $H^{2,0}(M,J)\in {\Bbb P}H^2(M, \C)$.
The map $\Per:\; \Teich \arrow {\Bbb P}H^2(M, \C)$ is 
called {\bf the period map}.
\ed

\Remark 
$\Per$ maps $\Teich$ into an open subset of a 
quadric, defined by
\[
\Perspace:= \{l\in {\Bbb P}H^2(M, \C)\ \ | \ \  q(l,l)=0, q(l, \bar l) >0\}.
\]
The manifold $\Perspace$ is called {\bf the period space} of $M$.
\er

As follows from Proposition \ref{_period_homo_Proposition_} below, 
$\Perspace = \frac{SO(b_2-3,3)}{SO(b_2-3,1)\times SO(2)}$.

\Theorem 
(Bogomolov, \cite{_Bogomolov_78_})\\
Let $M$ be a simple hyperk\"ahler manifold,
and $\Teich$ its Teichm\"uller space. Then
the period map $\Per:\; \Teich \arrow \Per$ is etale
(has invertible differential everywhere).
\et

\Remark Bogomolov's theorem implies that
$\Teich$ is smooth. It is non-Hausdorff
even in the simplest examples.
\er

Now the global Torelli theorem can be stated as follows.
Recall that the birational Teichm\"uller space $\Teich_b$ is a Hausdorff
reduction of the Teichm\"uller space of the holomorphic
symplectic manifolds of K\"ahler type.

\Theorem\label{_Glo_Tor_Theorem_}
Let $M$ be a simple hyperk\"ahler manifold,
and $\Per:\; \Teich_b \arrow \Perspace$ the period map.
Then $\Per$ is a diffeomorphism on each connected component.
\et

The following proposition is proven in a straghtforward
manner using 1950-ies style arguments of geometric topology.

\Proposition[(The Covering Criterion)]
Let $X\stackrel \phi \arrow Y$ be an etale map of smooth manifolds.
Suppose that each $y\in Y$ has a neighbourhood $B\ni y$ diffeomorphic
to a closed ball, such that  for each connected component
$B' \subset \phi^{-1}(B)$, $B'$ projects to $B$ surjectively.
Then $\phi$ is a covering.
\ep

Now, the Global Torelli implied by the following 
result, which is proven in Subsection \ref{_moduli_twistor_Subsection_}
using hyperk\"ahler structures.

\Proposition\label{_Covering_Proposition_}
 In assumptions of 
Theorem \ref{_Glo_Tor_Theorem_}, 
the period map satisfies the conditions of the
Covering Criterion.
\ep

\subsection{Moduli of hyperk\"ahler structures and twistor curves}
\label{_moduli_twistor_Subsection_}

\Proposition\label{_period_homo_Proposition_}
The period space 
\[
\Perspace:= \{l\in {\Bbb P}H^2(M, \C)\ \ | \ \  q(l,l)=0, q(l, \bar l) >0\}
\]
is identified with $\frac{SO(b_2-3,3)}{SO(2) \times SO(b_2-3,1)}$, which
is a Grassmannian of positive oriented 2-planes in $H^2(M,\R)$.
\ep

{\bf 
Proof. Step 1:} Given $l\in {\Bbb P}H^2(M, \C)$,  the space
generated by $\Im l, \Re l$ is 2-dimensional, because 
$q(l,l)=0, q(l, \bar l)$ implies that $l \cap H^2(M,\R)=0$.

{\bf  Step 2:}  This 2-dimensional plane is 
positive, because 
 $q(\Re l, \Re l) = q(l+ \bar l, l+ \bar l) = 2 q(l, \bar l)>0$.

{\bf Step 3:} Conversely, for any 2-dimensional positive
plane  $V\in H^2(M,\R)$, 
the quadric $\{l\in V \otimes_\R \C\ \ |\ \ q(l,l)=0\}$
consists of two lines; a choice of a line is determined by orientation.
\endproof

\Remark
Two hyperk\"ahler structures $(M,I,J,K,g)$
and $(M,I',J',K',g)$ are called {\bf equivalent}
if there exists a unitary quaternion $h$ such that
$I'=hIh^-1$,  $J'=hJh^-1$,  $K'=hKh^-1$. From the
holonomy characterization of simple hyperk\"ahler
manifolds (Remark \ref{_simple_holo_Remark_}) it
follows that two hyperk\"ahler structures are isometric
if and only if they are equivalent.
\er

\Definition
Let $(M,I,J,K,g)$ be a hyperk\"ahler manifold.
{\bf A hyperk\"ahler 3-plane} in $H^2(M,\R)$ is a positive oriented
3-dimensional subspace $W$, generated by three K\"ahler
forms $\omega_I, \omega_J, \omega_K$.
\ed

\Definition
Similarly to the Teichm\"uller space and period map
of complex structures, one can define the period
space of hyperk\"ahler metrics. Denote it by
$\Teich_H$. The corresponding period map is
\[ \Per:\; \Teich_H \arrow \Perspace_H,
\]
where $\Perspace_H=\frac{SO(b_2-3,3)}{SO(3) \times SO(b_2-3)}$
is the space of positive, oriented 3-planes, and
$\Per$ maps a hyperk\"ahler structure to the 
corresponding hyperk\"ahler 3-plane.
\ed

\Remark
There is one significant difference between 
$\Teich$ and the hyperk\"ahler Teichm\"uller space
$\Teich_H$: the latter is Hausdorff, and, in fact,
metrizable. Indeed, we could equip the space
 $\Teich_H$ of hyperk\"ahler metrics with the
Gromov-Hausdorff metric.
\er

Let $I\in \Teich$ be a complex structure, and ${\cal K}(I)$
its K\"ahler cone. The set of hyperk\"ahler metrics
compatible with $I$ is parametrized by ${\cal K}(I)$, by Calabi-Yau theorem.
 The corresponding 3-dimensional subspaces are 
generated by $\Per(I)+\omega$, where $\omega\in {\cal K}(I)$.
The local Torelli theorem implies that locally
$I\in \Teich$ is uniquely determined 
by the 2-plane generated by $\omega_J$ and $\omega_K$;
Calabi-Yau theorem implies that the hyperk\"ahler
metric is uniquely determined by the complex
structure and the K\"ahler structure. This 
gives the following hyperk\"ahler version of the
local Torelli theorem.

\Theorem \label{_Teich_H_local_Torelli_Theorem_}
Let $M$ be a simple hyperk\"ahler manifold,
and $\Teich_H$ its hyperk\"ahler Teichm\"uller space. Then
the period map $\Per:\; \Teich \arrow \Per_H$ mapping an
equivalence class of hyperk\"ahler structures to is its 3-plane
is etale (has invertible differential everywhere). \endproof
\et

\Remark
Let $W\subset H^2(M,\R)$ be a positive 3-dimensional plane.
The set $S_W\subset \Perspace$ 
of oriented 2-dimensional planes in $W$ is identified 
with $S^2 =\C P^1$. When $W$ is a hyperk\"ahler
3-plane,  $S_W$ is called {\bf the twistor family}
of a hyperk\"ahler structure. A point in the twistor family
corresponds to a complex structure $aI + bJ + cK \in {\Bbb H}$,
with $a^2+b^2+c^2=1$. We call the corresponding rational
curves $\C P^1\subset \Teich$ {\bf the twistor lines}.
It is not hard to see that the twistor lines are holomorphic.
\er

\Definition
Let $W\in \Perspace_H$ be a positive 3-plane, 
$S_W\subset \Perspace$ the corresponding rational curve,
and $x\in S_W$ be a point.
It is called {\bf liftable} if for any point 
$y\in \Per^{-1}(x)\subset \Teich$
there exists ${\cal H}\in \Teich_H$
such that the corresponding twistor line
contains $y$.
\ed

When $W$ is generic, the corresponding 
line $S_W$ is liftable, as indicated below.

\newcommand{\NS}{\operatorname{\sf NS}}

\Definition
{\bf The Neron-Severi lattice} $NS(I)$ of a hyperk\"ahler manifold
$(M,I)$ is $H^{1,1}(M,I)\cap H^2(M,\Z)$.
\ed

The following theorem, based on 
results of \cite{_Demailly_Paun_}, was proven by D. Huybrechts.

\Theorem[(\cite{_Huybrechts:cone_})]
\label{_gene_cone_Theorem_}
Let $M$ be a hyperkaehler manifold with 
$\NS(M)=0$. Then its Kaehler cone is
one of two components of the set 
\[ \{ \nu \in H^{1,1}(M,\R) \ | \ q(\nu, \nu)\geq 0\}.\]
\et

\Definition
Let $S\subset \Teich$ be a $\C P^1$ 
associated with a twistor family. It is called {\bf
generic} if it passes through a point $I\in \Teich$
with $\NS(M,I)=0$. Clearly,
a hyperk\"ahler 3-plane $W\subset H^2(M,\R)$
corresponds to a generic twistor family if and only if
its orthogonal complement $W^\bot\subset H^2(M,\R)$
does not contain rational vectors. 
A 3-plane $W\in \Perspace_H$ is called {\bf generic}
if $W^\bot\subset H^2(M,\R)$ does 
not contain rational vectors. The corresponding
rational curve $S_W\subset \Perspace$ is called
{\bf a GHK line}. GHK lines are liftable, which
is very useful for many purposes, including
the proof of Torelli theorem (see also 
\cite{_Amerik_Verbitsky:rational_curv_}, 
where GHK lines were used to study K\"ahler 
cones of hyperk\"ahler manifolds).
\ed

The following theorem immediately follows
from the Calabi-Yau theorem and the description
of the K\"ahler cone given in Theorem \ref{_gene_cone_Theorem_}.

\Theorem\label{_liftable_Theorem_}
Let $W\in \Perspace_H$ be a generic plane,
$S_W\subset \Perspace$ the corresponding 
rational curve, and $x\in S_W$ a generic point.
Then $(S_W,x)$ is liftable. \endproof
\et

Assumptions of the
covering criterion (Proposition \ref{_Covering_Proposition_})
immediately follow from Theorem \ref{_liftable_Theorem_}.
Indeed, it is not hard to see that 
any two points on a closed ball $B\subset \Perspace$
can be connected inside $B$ by a sequence of GHK curves intersecting
in generic points of $B$. Since these curves are liftable, 
any connected component of $\Per^{-1}(B)$ is mapped
to $B$ surjectively.


\section{Teichm\"uller spaces and ergodic theory}


\subsection{Ergodic complex structures}

After the Teichm\"uller space and the mapping class group 
are understood, it is natural to 
consider the quotient space $\Comp/\Diff=\Teich/\Gamma$
of the Teichm\"uller space by the mapping class group
$\Gamma:=\Diff/\Diff_0$. 

\Claim
Let $M$ be a simple hyperk\"ahler manifold,
$\Gamma$ its mapping class group, and $\Teich_b$ the
birational Teichm\"uller space. Then the quotient
$\Teich_b/\Gamma$ parametrizes the birational
classes of deformations of $M$.
\ec

One could call the quotient $\Teich_b/\Gamma$
``the moduli space'', but, unfortunately, this
is not a space in any reasonable sense. Indeed, as we shall
see, non-trivial closed subsets of $\Teich_b/\Gamma$
are at most countable, making $\Teich_b/\Gamma$ terribly
non-Hausdorff. This means that the concept of ``moduli space''
has no meaning, and all interesting information about
moduli problems is hidden in dynamics of $\Gamma$-action 
on $\Teich$.

Let $I\in \Teich$ be a point, and $\Teich^I\subset \Teich$ its connected
component. Since $\Teich$ has finitely many components,
a subgroup mapping class group fixing $\Teich$ has finite index.
Its image in $\Aut(\Teich^I)$ is called {\bf monodromy group}
and denoted $\Gamma^I$ (Definition \ref{_Monodro_Definition_}).
It is a finite index subgroup in $SO(H^2(M, \Z))$.

All that said, we find that the moduli problem for hyperk\"ahler manifold
is essentially reduced to the dynamics of the $\Gamma^I$-action on the
space $\Perspace$, which is understood as a 
Grassmannian of positive, oriented 2-planes in $H^2(M,\R)$.

It is natural to study the dynamics of a group action from the
point of view of ergodic theory, ignoring measure zero subsets.
However, the quotient map $\Teich\arrow \Teich_b$
is bijective outside of a union of countably many divisors,
corresponding to complex structures $I$ with $NS(M,I)$ non-zero.
This set has measure 0. Therefore, the quotient map
$\Teich\arrow \Teich_b$ induces an equivalence of measured spaces.
For the purposes of ergodic theory, we shall identify 
$\Teich^I$ with the corresponding homogeneous space
$\Perspace$.

Ths first observation, based on a theorem of C. Moore, implies that 
the monodromy action on $\Perspace$ is ergodic.

\Definition
Let $(M,\mu)$ be a space with measure,
and $G$ a group acting on $M$.
This action is {\bf ergodic} if all
$G$-invariant measurable subsets $M'\subset M$
satisfy $\mu(M')=0$ or $\mu(M\backslash M')=0$.
\ed

\Claim
Let $M$ be a manifold, $\mu$ a Lebesgue measure, and
$G$ a group acting on $(M,\mu)$ ergodically. Then the 
set of non-dense orbits has measure 0.
\ec

{\bf Proof:} Consider a non-empty open subset $U\subset M$. 
Then $\mu(U)>0$, hence $M':=G\cdot U$ satisfies 
$\mu(M\backslash M')=0$. For any orbit $G\cdot x$
not intersecting $U$, one has $x\in M\backslash M'$.
Therefore, the set of such orbits has measure 0.
\endproof

\Definition
Let $I\in \Comp$ be a complex structure on a manifold.
It is called {\bf ergodic} if its $\Diff$-orbit
is dense in its connected component of $\Comp$.
\ed

\Remark
This is equivalent to density of $\Gamma$-orbit
of $I$ in its Teichm\"uller component.
\er

\subsection{Ergodicity of the monodromy group action}
\label{_ergodi_theorem_Subsection_}

\Definition 
Let $G$ be a Lie group, and $\Gamma\subset G$ a discrete
subgroup. Consider the pushforward of the Haar measure to
$G/\Gamma$. We say that $\Gamma$ {\bf has finite covolume}
if the Haar measure of $G/\Gamma$ is finite.
In this case $\Gamma$ is called {\bf a lattice subgroup}.
\ed

\Remark
Borel and Harish-Chandra proved that
an arithmetic subgroup of a reductive group $G$
is a lattice whenever $G$ has no non-trivial characters
over $\Q$ (see e.g. \cite{_Vinberg_Gorbatsevich_Shvartsman_}). 
In particular, all arithmetic subgroups
of a semi-simple group are lattices.
\er

\Theorem (Calvin C. Moore, \cite[Theorem 7]{_Moore:ergodi_})\\
Let $\Gamma$ be an arithmetic subgroup in a non-compact 
simple Lie group $G$ with finite center, and $H\subset G$ a 
non-compact subgroup. Then the left action of $\Gamma$
on $G/H$ is ergodic.
\et

\Theorem Let ${\Teich}$ be a connected component of 
a Teichm\"uller space, and
$\Gamma^I$ its monodromy group. Then the set of all
non-ergodic points of $\Teich$ has measure 0.
\et

{\bf  Proof:} Global Torelli theorem 
identifies $\Teich$ (as a measured space) and $G/H$, where
$G=SO(b_2-3,3)$, $H=SO(2) \times SO(b_2-3,1)$.
Since $\Gamma^I$ is an arithmetic lattice, 
$\Gamma^I$-action on $G/H$ is ergodic,
by Moore's theorem.
\endproof

\hfill

Moore's theorem implies that outside of a measure zero
set, all complex structures on $\Teich$ are ergodic.
If we want to determine which exactly complex
structures are ergodic, we have to use Ratner's theorem,
giving precise description of a closure of a $\Gamma^I$-orbit
in a homogeneous space.

Now I will state some basic results of Ratner theory.
For more details, please see \cite{_Kleinbock_etc:Handbook_} 
and \cite{_Morris:Ratner_}.

\hfill

\Definition
Let $G$ be a Lie group, and $g\in G$ any element.
We say that $g$ is {\bf unipotent} if $g=e^h$ for a
nilpotent element $h$ in its Lie algebra.
A group $G$ is {\bf generated by unipotents}
if $G$ is multiplicatively generated by unipotent elements.
\ed

\Theorem (\cite[1.1.15 (2)]{_Morris:Ratner_})\\
Let $H\subset G$ be a Lie subroup generated by 
unipotents, and $\Gamma\subset G$ a lattice.
Then a closure of any $H$-orbit in $G/\Gamma$
is an orbit of a closed, connected subgroup $S\subset G$,
such that $S\cap \Gamma\subset S$ is a lattice.
\et

When this lattice is arithmetic, one could describe
the group $S$ very explicitly.

\Claim ( \cite[Proposition 3.3.7]{_Kleinbock_etc:Handbook_}
or \cite[Proposition 3.2]{_Shah:uniformly_})\\
Let $x\in G/H$ be a point in a homogeneous space,
and $\Gamma\cdot x$ its $\Gamma$-orbit, where $\Gamma$
is an arithmetic lattice. Then
its closure is an orbit of a group $S$ containing
stabilizer of $x$. Moreover, $S$ is a smallest
group defined over rationals and stabilizing $x$.
\ec

For the present purposes, we are interested
in a pair $SO(3,k)\supset SO(1,k)\times SO(2)\subset G$
(or, rather, their connected components $G=SO^+(3,k)$
and $H= SO(1,k)\times SO(2)\subset G$).
In this case, there are no intermediate subgroups.

\Claim Let $G=SO^+(3,k)$, and 
$H\cong SO^+(1,k)\times SO(2)\subset G$.
Then any closed connected
Lie subgroup $S\subset G$ containing $H$ coincides
with $G$ or with $H$.
\ec

\Corollary
Let $J\in \Perspace=G/H$. Then either $J$ is ergodic, or
its $\Gamma$-orbit is closed in $\Perspace$. 
\endproof
\eco

By Ratner's theorem, in the latter case
the $H$-orbit of $J$ has finite volume in $G/\Gamma$.
Therefore,  its intersection with $\Gamma$ is a lattice in 
$H$. This brings

\Corollary
Let $J\in \Perspace$ be a point such that its 
$\Gamma$-orbit is closed in $\Perspace$. Consider its stabilizer
$\St(J)\cong H \subset G$. Then $\St(J)\cap \Gamma$
is a lattice in $\St(J)$.
\endproof
\eco

\Corollary \label{_ergodic_Pic_max_Corollary_}
Let $J$ be a non-ergodic complex structure on a hyperk\"ahler
manifold, and $W\subset H^2(M,\R)$ be a plane generated
by $\Re\Omega, \Im \Omega$. Then $W$ is rational.
Equivalently, this means that $Pic(M)$ has maximal possible
dimension.
\endproof
\eco

Similar results are true for a torus of dimension $>1$;
we refer the reader to \cite{_V:ergodic_} 
for precise statements and details of the proof.


\section{Applications of ergodicity}


\subsection{Ergodic complex structures, 
Gromov-Hausdorff closures, and semicontinuity}

The ergodicity theorem (Subsection
\ref{_ergodi_theorem_Subsection_})
has some striking and even paradoxical implications.
For instance, consider a K\"ahler cone $\Kah$ of a
hyperk\"ahler manifold (or a torus of dimension $>1$)
equipped with an ergodic complex structure.
By Calabi-Yau theorem, each point of $\Kah$ corresponds to
a Ricci-flat metric on $M$. If we restrict ourselves to
those metrics which satisfy $\diam(M,g)\leq d$
(with bounded diameter), then, by Gromov's compactness
theorem (\cite{_Gromov:Riemannian_}), the set $X_d$
of such metrics is precompact in the Gromov's space
of all metric spaces, equipped with the Gromov-Hausdorff
metric. It is instructive to see what kind of metric
spaces occur on its boundary (that is, on $\bar
X_d\backslash X_d$). To see this, let $\nu_i$ be a sequence
of diffeomorphisms satisfying $\lim_i\nu_i(I)=I'$.
By Kodaira stability theorem, the K\"ahler cone of $(M,I)$ is lower continuous
on $I$. Therefore, there exists a family of K\"ahler classes
$\omega_i$ on $(M, \nu_i(I))$ which converge to a given
K\"ahler class $\omega'$ on $(M,I')$. This implies
convergence of the corresponding Ricci-flat metrics.
We obtain that any Ricci-flat metric on $(M,I')$
(for any $I'$ in the same deformation class as $I$)
is obtained as a limit of Ricci-flat metrics on $(M,I)$.

This gives the following truly bizzarre theorem.

\Theorem
Let $(M,I)$ be an ergodic complex structure on a
hyperk\"ahler manifold, $X\cong\Kah$ the set of all
Ricci-flat K\"ahler metrics on $(M,I)$, 
and $g'$ another Ricci-flat metric
on $M$ in the same deformation class. Then $g'$
lies in the closure of $X$ with respect to the
Gromov topology on the space of all metrics.
\et

This result is very strange, because $\Kah$ 
is a smooth manifold of dimension $b_2(M)-2$. 
By Theorem \ref{_Teich_H_local_Torelli_Theorem_}, 
the space $\Teich_H$ of all hyperk\"ahler
metrics is a smooth manifold of 
dimension $\frac{b_2(b_2-1)(b_2-2)}6$,
clearly much bigger than $\dim \Kah$.
Obviously, the boundary of $X$ is highly 
irregular and chaotic.

For another application, consider some numerical quantity 
$\mu(I)$ associated with an equivalence class of complex
structures. Suppose that $\mu$ is continuous or
semi-continuous on $\Teich$. Then $\mu$ is constant 
on ergodic complex structures. To see this,
suppose that $\mu$ is upper semicontinuous,
giving 
\begin{equation}\label{_semiconti_Equation_}
\mu(\lim_k I_k)\geq \lim_k(\mu(I_k)).
\end{equation}
Given an ergodic complex structire $I$,
find a sequence $I_k=\nu_k(I)$ converging
to a complex structure $I'$.
Then \eqref{_semiconti_Equation_} gives $\mu(I)\leq
\mu(I')$.
This implies that any ergodic complex structure 
satisfies $\mu(I)=\inf_{I'\in \Teich} \mu(I')$.

This observation can be applied to Kobayashi pseudometric
and Kobayashi hyperbolicity.

\subsection{Kobayashi non-hyperbolicity of hyperk\"ahler manifolds}

\Definition
{\bf Pseudometric} on $M$ 
is a function $d:\; M \times M \arrow \R^{\geq 0}$
which is symmetric: $d(x,y)=d(y,x)$ and satisfies the
triangle inequality $d(x,y)+d(y,z) \geq d(x,z)$.
\ed

\Remark 
Let ${\goth D}$ be a set of pseudometrics.  Then
$d_{\max}(x,y):= \sup_{d\in {\goth D}}d(x,y)$ is also a
pseudometric.
\er

\Definition
The  Kobayashi pseudometric on a complex manifold $M$
is $d_\max$ for the set ${\goth D}$ of all pseudometrics
such that any holomorphic map from the 
Poincar\'e disk to $M$ is distance-non-increasing.
\ed

In other words, a Kobayashi pseudo-distance
between two points $x,y$ is an infimum of distance
from $x$ to $y$ in Poincare metric for any sequence of holomorphic
disks connecting $x$ to $y$.

The following observation is not difficult to see.

\Claim
Let $\pi:\; {\cal M} \arrow X$ be a smooth holomorphic family,
which is trivialized as a smooth manifold: ${\cal M}=M \times X$,
and $d_x$ the Kobayashi metric on $\pi^{-1}(x)$. Then $d_x(m,m')$
is upper continuous on $x$. \endproof
\ec

\Corollary
Denote the diameter of the Kobayashi pseudometric by
$\diam(d_x):= \sup_{m,m'}d_x(m,m')$.  
Then the Kobayashi diameter of a fiber of $\pi$ 
is an upper continuous function: 
$\diam:\; X \arrow \R^{\geq 0}$. \endproof
\eco

For a projective K3 surface, the Kobayashi pseudometric
vanishes (\cite[Lemma 1.51]{_Voisin:kobayashi_}).
However, all non-projective K3 surfaces are ergodic
(Corollary \ref{_ergodic_Pic_max_Corollary_}).
This proves the vanishing of Kobayashi pseudodistance
for all K3 surfaces. A more general version of this
result is due to 
due to Kamenova-Lu-Verbitsky.

\Theorem[(\cite{_Kamenova_Lu_Verbitsky_})]
Let $M$ be a Hilbert scheme of K3. Then the Kobayashi
pseudometric on $M$ vanishes
\et

\Definition
A complex manifold is called {\bf Kobayashi hyperbolic}
if the Kobayashi pseudometric is a metric.
\ed

\Definition
An entire curve is a non-constant 
map $\C \arrow M$.
\ed

Brody has shown that a compact manifold is Kobayashi
hyperbolic if and only if it admits no entire curves.
The same argument also proves semicontinuity.

\Theorem (\cite{_Brody:hyperbolic_})
Let $I_i$ be a sequence of complex structures on $M$
which are not hyperbolic, and $I$ its limit. Then
$(M,I)$ is also not hyperbolic.
\et

With ergodicity, this can be used to prove that all hyperk\"ahler
manifolds are non-hyperbolic.

Recall that {\bf a twistor family} of complex structures
on a hyperk\"ahler manifold $(M,I,J,K)$ is a family 
of complex structures of form 
$S^2 \cong \{ L:= aI + bJ +c K, \ \ \ a^2+b^2+c^2=1\}$.
F. Campana has obtained a remarkable partial result
towards non-hyperbolicity.

\Theorem (\cite{_Campana:twistor_hy_})
Let $M$ be a hyperk\"ahler manifold, and 
$S\subset \Teich$ a twistor family.
Then there exists an entire curve in some
$I\in S$.
\et

\Claim\label{_only_ergo_twi_Claim_}
There exists a twistor family which has only ergodic
fibers. 
\ec

{\bf Proof:} There are only countably many complex
structures which are not ergodic; however, twistor curves
move freely through the Teichm\"uller space of a hyperk\"ahler manifold,
as seen from Theorem \ref{_liftable_Theorem_}. \endproof

\hfill

Applying Campana's theorem to the family
constructed in Claim \ref{_only_ergo_twi_Claim_}, we obtain
an ergodic complex structure which is non-hyperbolic.
Then the Brody's theorem implies that all complex structures
in the same deformation class are non-hyperbolic.

\Theorem
All hyperk\"ahler manifolds are non-hyperbolic.
\endproof
\et

\subsection{Symplectic packing and ergodicity}

I will finish this talk with a list of open problems
of hyperk\"ahler and holomorphically symplectic geometry
which might be solvable with ergodic methods.

\Question
Let $M$ be a hyperk\"ahler manifold, and $\Teich$ its 
Teichm\"uller space.
Consider the universal fibration ${\cal X} \arrow \Teich$
(Remark \ref{_Universal_fi_Remark_}). The mapping class group
$\Gamma$ acts on ${\cal X}$ in a natural way. 
Is this action ergodic?
\eq

This question (suggested by Claire Voisin)
seems to be related to the following conjecture.

\Conjecture
Let $M$ be a K3 surface. Then for each $x\in M$ and
$v\in T_x M$ there exists an entire curve $C\ni x$
with $T_x C \ni v$.
\econ

The symplectic packing problem is a classical subject of
symplectic geometry (\cite{_McDuff-Polterovich:packing_}). 
However, its holomorphically symplectic version seems
to be completely unexplored.

\Definition
{\bf A holomorphic symplectic ball} $B_r$ of radius $r$
is a complex holomorphically symplectic manifold admitting
a holomorphic symplectomorphism to an open ball in $\C^{2n}$
of radius $r$ with the standard holomorphic symplectic
form $\sum_{i=1}^{n} dz_{2i-1}\wedge dz_{2i}$.
\ed

Notice that by a holomorphic symplectic version
of Darboux theorem, any holomorphically symplectic
manifold is locally symplectomorphic to a 
holomorphic symplectic ball.

\Definition
Let $M$ be a holomorphically symplectic manifold.
{\bf Symplectic packing} of radii $r_1, ..., r_k$ 
of $M$ is a set of holomorphic symplectomorphisms 
$\phi_i:\; B_{r_i}\hookrightarrow M$ with images
of $\phi_i$ not intersecting.
\ed

Obviously, in these assumptions, $\sum \Vol(B_{r_i})\leq \Vol_M$,
where $\Vol$ denotes the symplectic volume of a holomorphic
symplectic manifold $(M, \Omega_M)$:
\[ \Vol(M)=\int_M (\Omega_M\wedge\bar\Omega_M), \ \ 2n =\dim_\C M.
\]
The volume inequality puts certain restrictions on the
possible symplectic packing. Are there any other restrictions?

For the usual (smooth) symplectic packing, some
additional restrictions are obtained from the Gromov's
symplectic capacity theorem and from the study
of pseudoholomorphic curves. However, it seems
that in holomorphic symplectic situation these
restrictions are also trivial. For a general 
compact torus of real dimension 4, volume is known
to be the only restriction to existence of
symplectic packing (\cite{_McDuff_etc:torus_packings_}).
It seems that a similar result about the smooth symplectic
packings is true for K3 surfaces
as well, and, possibly, for any hyperk\"ahler manifold.

The arguments used to treat the usual (smooth)
symplectic packings don't work for the holomorphic
symplectic case. However, the set of possible radii
for symplectic packing is obviously semicontinous,
hence it can be studied by ergodic methods, in the
same way as one studies the Kobayashi pseudometric.

The following classical question was treated Buzzard and Lu
in \cite{_Buzzard_Lu:dominable_}.

\Definition
A complex manifold $M$ of dimension $n$ is called
{\bf dominated by $\C^n$} if there exists
a holomorphic map $\phi:\; \C^n \arrow M$
which has non-degenerate differential in generic point.
\ed

Buzzard and Lu proved that Kummer K3 surfaces are dominated
by $\C^2$. So far, there is not a single example of a hyperk\"ahler
manifold $M$ for which it is proven that $M$ is not dominated.
This leads to the following conjecture

\Conjecture
Any compact hyperk\"ahler manifold is dominated by $\C^n$.
\econ

There is no semicontinuity in dominance, because Brody lemma
fails to produce dominating maps $\C^n \arrow M$ for $n>1$ as limits
of sequences of dominating maps. In the proof of Brody's lemma
(showing that a limit of a sequence of 
entire curves contains an entire curve)
one takes a reparametrizations of each of the curves
in the sequence. Starting from a sequence of
dominating maps, one could apply the same argument,
but each subsequent reparametrization can lead to
smaller Jacobian of the differential, and the differential
of the limit could be zero.

It seems that more of the Brody's argument can be retained
if we restrict ourselves to symplectomorphisms.

\Question
Consider a flat holomorphically symplectic structure on $\C^2$.
Is there a holomorphic map $\C^2 \arrow M$ to a K3
surface which is compatible with the holomorhic symplectic form?
\eq

Probably not. However, a quantitative version of this question
makes sense. Let $M$ be a hyperk\"ahler manifold,
and $K(M)$ the supremum of all $r$ such that
there exists a symplectic immersion from a
symplectic ball of radius $r$ to $M$. It is not hard
to see that $K(M)$ is semicontinuous in families,
hence constant on ergodic complex structures.

\Question
For a given hyperk\"ahler manifold,
find $K(M)$.
\eq

It is not clear if $K(M)$ is finite or infinite, even for a K3 surface
(it is clearly infinite for a torus).



\begin{thebibliography}{99999}


\bibitem[AV]{_Amerik_Verbitsky:rational_curv_}
Ekaterina Amerik, Misha Verbitsky
{\em Rational curves on hyperkahler manifolds,}
 arXiv:1401.0479, 34 pages. 

\bibitem[Bea1]{_Beauville:hk_} 
 Beauville, A. {\em 
Varietes K\"ahleriennes dont la premi\`ere classe de Chern est
nulle.}  J. Diff. Geom. {\bf 18}, pp. 755-782 (1983).


\bibitem[Bea2]{_Beauville:K3_}
Beauville, Arnaud
{\em Le theoreme de Torelli pour les surfaces $K3$: fin de la demonstration}
Geometry of $K3$ surfaces: moduli and periods (Palaiseau, 1981/1982).
Asterisque No. 126 (1985), 111-121. 

\bibitem[Bes]{_Besse:Einst_Manifo_} 
Besse, 
A., {\em Einstein Manifolds}, Springer-Verlag, New York (1987)



\bibitem[B:1974]{_Bogomolov:decompo_}  
Bogomolov, F. A., {\em On the decomposition of 
K\"ahler manifolds with trivial canonical class}, Math. USSR-Sb.
{\bf 22} (1974), 580-583.


\bibitem[B:1978]{_Bogomolov_78_}
Bogomolov, F. A. {\bf \it Hamiltonian K\"ahlerian manifolds. }
 Dokl. Akad. Nauk SSSR 243 (1978), no. 5, 1101-1104.



\bibitem[B:1981]{_Bogomolov_81_} Bogomolov, F. A.,
{\em K\"ahler 
manifolds with trivial canonical class,} Preprint, Institut
des Hautes Etudes Scientifiques (1981), 1-32.


\bibitem[Br:1978]{_Brody:hyperbolic_}
 Brody, R., {\em Compact manifolds and hyperbolicity,}
Trans. Amer. Math. Soc. 235 (1978), 213-219.


 \bibitem[BL]{_Buzzard_Lu:dominable_}
        Buzzard, G.; Lu, S. S.-Y.,
        {\em Algebraic surfaces holomorphically dominable by $\C^2$,}
        to appear in Invent. Math.


\bibitem[Cam]{_Campana:twistor_hy_}
 F. Campana,
{\em An application of twistor theory to 
the nonhyperbolicity of certain compact symplectic 
K\"ahler manifolds} 
J. Reine Angew. Math., 425:1-7, 1992.


\bibitem[C]{_Catanese:moduli_}
F. Catanese, 
{\em A Superficial Working Guide to Deformations and Moduli},
arXiv:1106.1368, 56 pages.



\bibitem[De]{_Debarre:Torelli_}
Debarre, O.,
{\em Un contre-exemple au th\'eor\`eme 
de Torelli pour les vari\'et\'es symplectiques irr\'eductibles,}
C. R. Acad. Sci. Paris S\'er. I Math. 299 (1984), no. 14, 681--684. 



\bibitem[DP]{_Demailly_Paun_}
Demailly, J.-P., Paun, M., 
{\em Numerical characterization of the K\"ahler cone 
of a compact K\"ahler manifold}, 
Annals of Mathematics, 
{\bf 159} (2004), 1247-1274, math.AG/0105176.


\bibitem[Dou]{_Douady:Bourbaki_}
Douady, A., {\em 
Le probleme des modules pour les varietes analytiques complexes,}
   Seminaire Bourbaki, 1964/1965, No 277.

\bibitem[Fr]{_Friedman_}
Friedman, Robert,
{\em A new proof of the global Torelli theorem for $K3$ surfaces},
Ann. of Math. (2) 120 (1984), no. 2, 237-269. 


\bibitem[Fu]{_Fujiki:HK_}  
Fujiki, A. {\em On the de Rham Cohomology Group of a Compact 
K\"ahler Symplectic Manifold}, Adv. Stud.
Pure Math. 10 (1987), 105-165.




\bibitem[Gr]{_Gromov:Riemannian_} 
Gromov, Misha, {\em  Metric structures for Riemannian and
non-Riemannian spaces}, Based on the 1981 French
original. With appendices by M. Katz, P. Pansu and
S. Semmes. Translated from the French by Sean Michael
Bates. Progress in Mathematics, 152. Birkh\"auser Boston,
Inc., Boston, MA, 1999. xx+585 pp.



\bibitem[Ha]{_Hamilton:inverse_}
Hamilton, Richard S. {\em The inverse function theorem of Nash and Moser,} 
Bull. Amer. Math. Soc. (N.S.) 7 (1982), no. 1, 65-222. 


\bibitem[Hu1]{_Huybrechts:cone_}
Huybrechts, D.
{\em The K\"ahler cone of a compact hyperk\"ahler manifold},
Math. Ann. 326 (2003), no. 3, 499--513, arXiv:math/9909109.



\bibitem[Hu2]{_Huybrechts:finiteness_}
Huybrechts, D., 
{\em Finiteness results for hyperk\"ahler manifolds},
 J. Reine Angew. Math.  558  (2003), 15--22, arXiv:math/0109024.

\bibitem[KLV]{_Kamenova_Lu_Verbitsky_}
Ljudmila Kamenova, Steven Lu, Misha Verbitsky,
{\em Kobayashi pseudometric on hyperkahler manifolds},
arXiv:1308.5667, 21 pages.

\bibitem[KSS]{_Kleinbock_etc:Handbook_}
Kleinbock, Dmitry; Shah, Nimish; Starkov, Alexander,
{\em Dynamics of subgroup actions on homogeneous spaces of
  Lie groups and applications to number theory}, Handbook
of dynamical systems, Vol. 1A, 813-930, North-Holland,
Amsterdam, 2002. 


\bibitem[Kul]{_Kulikov_}
Kulikov,  Vik. S.,
{\em 
Degenerations of K3 surfaces and Enriques surfaces,} 
Mathematics of the USSR-Izvestiya, 11:5 (1977), 957-989.


\bibitem[Ku1]{_Kuranishi:new_}
Kuranishi, M.
{\em New proof for the existence of locally complete 
families of complex structures}, 1965 Proc. Conf. 
Complex Analysis (Minneapolis, 1964) pp. 142-154 Springer, Berlin.

\bibitem[Ku2]{_Kuranishi:note_}
Kuranishi, Masatake
{\em A note on families of complex structures} 1969 Global
Analysis (Papers in Honor of K. Kodaira) pp. 309-313
Univ. Tokyo Press, Tokyo.

\bibitem[LMS]{_McDuff_etc:torus_packings_}
Janko Latschev, Dusa McDuff, Felix Schlenk, 
{\bf The Gromov width of 4-dimensional tori},
 	arXiv:1111.6566.


\bibitem[L]{_Loojenga_}
E. Loojenga,
{\em A Torelli theorem for K\"ahler-Einstein K3 surfaces},
Lecture Notes in Mathematics, 1981, Volume 894/1981, 107-112.

\bibitem[MP]{_McDuff-Polterovich:packing_} D. McDuff and L. Polterovich, 
{\em Symplectic packings and algebraic geometry,} Invent.
     Math. 115 (1994), 405-434.



\bibitem[Ma1]{_Markman:mono_}
Markman, E., {\em On the monodromy of moduli spaces of sheaves
on $K3$ surfaces},  J. Algebraic Geom.  17  (2008),  no. 1,
29--99, arXiv:math/0305042.

\bibitem[Ma2]{_Markman:constra_}
Markman, E. {\em
Integral constraints on the monodromy group of
    the hyperkahler resolution of a symmetric product of a
    K3 surface,} International Journal of Mathematics
Vol. 21, No. 2 (2010) 169-223, arXiv:math/0601304.
 



\bibitem[Mo:1966]{_Moore:ergodi_}
Calvin C. Moore,
{\em Ergodicity of Flows on Homogeneous Spaces},
American Journal of Mathematics
Vol. 88, No. 1 (Jan., 1966), pp. 154-178


\bibitem[Mor]{_Morris:Ratner_}
Morris, Dave Witte, 
{\em Ratner's Theorems on Unipotent Flows,} 
Chicago Lectures in Mathematics, University 
of Chicago Press, 2005.

\bibitem[Na]{_Namikawa:Torelli_}
Namikawa, Y.,
{\em Counter-example to global Torelli problem for 
irreducible symplectic manifolds}, 
 Math. Ann.  324  (2002),  no. 4, 841--845.



\bibitem[O]{_O_Grady_}
O'Grady, Kieran G., 
{\em A new six-dimensional irreducible symplectic variety},
J. Algebraic Geom. 12 (2003), no. 3, 435--505.


\bibitem[PS]{_Piatetski_Shapiro_Shafarevich_}
Piatecki-Shapiro, I.I.; Shafarevich I.R.,
{\em Torelli's
theorem for algebraic surfaces of type
K3}, Izv. Akad. Nauk SSSR Ser. Mat. (1971) 35: 530-572.



\bibitem[Sh]{_Shah:uniformly_}
N. A. Shah,
{\em Uniformly distributed orbits of certain flows on homogeneous spaces,} 
Math. Ann. 289 (2) (1991), 315-33.



\bibitem[Su]{_Sullivan:infinite_}
Sullivan, D.,
{\em Infinitesimal computations in topology}, Publications
Math\-\'ema\-tiques de l'IH\'ES, 47 (1977), p. 269-331

\bibitem[Te:1944]{_Teichmuller:1944_}
Teichm\"uller, Oswald,
{\em Ver\"anderliche Riemannsche Fl\"achen}, 
Abh. Preuss. Deutsche Mathematik 7: 344-359,
1944.


\bibitem[Ti]{_Tian_} 
G. Tian, {\em Smoothness of the universal
deformation space of compact Calabi-Yau manifolds and its
Petersson-Weil metric}, in {\em Math. Aspects of String Theory},
S.-T. Yau, ed., Worlds Scientific, 1987, 629--646.



\bibitem[Tj]{_Tjurina_}
 Tjurina, G. N., 
{\em The space of moduli of a complex surface with q = 0
and K = 0}. In: "Algebraic Surfaces", Seminar Shafarevich, 
Proc. Steklov Inst. 75 (1965).



\bibitem[Tod1]{_Todorov:K3_}
Todorov, A. N. 
{\em Applications of the K\"ahler-Einstein-Calabi-Yau metric to 
moduli of K3 surfaces}, Inventiones Math. 6-1, 251-265 (1980).

\bibitem[Tod2]{_Todorov_} 
A. Todorov, {\em The Weil-Petersson geometry of
the moduli space of $SU(n \geq 3)$ (Calabi-Yau) manifolds}, Comm.,
Math. Phys. {\bf 126} (1989), 325--346.

\bibitem[To:1913]{_Torelli_}
Ruggiero Torelli, 
{\em Sulle variet\`a di Jacobi}, 
Rend. della R. Acc. Nazionale dei Lincei , (5), 22, 1913, 98-103.

\bibitem[Si]{_Siu:K3_}
Siu, Y.-T., {\em Every K3 surface is K\"ahler}, Invent. Math. 73
(1983), 139-150.



\bibitem[V:1996]{_Verbitsky:coho_announce_} 
Verbitsky, M.,
{\it Cohomology of compact hyperk\"ahler manifolds
and its applications,} GAFA vol. 6 (4) pp. 601-612 (1996).


\bibitem[V:2013]{_V:Torelli_}
Verbitsky, M.,
{\em A global Torelli theorem for hyperk\"ahler manifolds,}
 Duke Math. J. Volume 162, Number 15 (2013), 2929-2986.

\bibitem[V]{_V:ergodic_}
Verbitsky, M.,
{\em Ergodic complex structures on hyperkahler manifolds},
arXiv:1306.1498, 22 pages.



\bibitem[VGS]{_Vinberg_Gorbatsevich_Shvartsman_}
Vinberg, E. B.,  Gorbatsevich, V. V.,  Shvartsman, O. V., 
{\em Discrete Subgroups of Lie Groups}, in 
``Lie Groups and Lie Algebras II'',
Springer-Verlag, 2000.



\bibitem[Vo]{_Voisin:kobayashi_}
Claire Voisin,
{\em On some problems of Kobayashi and Lang; algebraic
approaches,}  Current Developments in Mathematics 2003,
no. 1 (2003), 53-125.


\bibitem[W:1957]{_Weil_} 
A. Weil. {\em Zum Beweis des Torellischen
Satzes,} Nachr. Akad. Wiss. G\"ottingen,
Math.-Phys. Kl. IIa: 32-53 (1957).


\end{thebibliography}
\end{document}